\documentclass[12pt,reqno]{amsart}

\usepackage{amsthm}
\usepackage{pstricks,pst-node,multido}
\usepackage{amssymb}
\usepackage{scrtime}

\DeclareMathOperator{\Aut}{Aut}
\DeclareMathOperator{\GL}{GL}

\DeclareMathOperator{\shape}{shape}

\newcommand{\PG}{PG}

\usepackage[left=35mm, right=35mm, top=25mm, bottom=20mm, a4paper, includefoot]{geometry}

\newcounter{daana} \newcounter{daanb}

\newcommand{\lista}[2]
{\begin{list}{#1}{%
\usecounter{daana}%
\setlength{\listparindent}{\parindent}%
\setlength{\leftmargin}{3.5em}%
\setlength{\itemindent}{0mm}%
\setlength{\labelsep}{.6em}%
\setlength{\labelwidth}{10em}%
\setlength{\itemsep}{0em}%
\setlength{\parsep}{.5mm}%
\setlength{\topsep}{1.5ex}%
\setlength{\listparindent}{\parindent}%
#2}}

\newcommand{\listb}[2]
{\begin{list}{#1}{%
\usecounter{daanb}%
\setlength{\rightmargin}{\leftmargin}%
\setlength{\listparindent}{\parindent}%
#2}}

\newcommand{\thspace}{3mm}

\newtheoremstyle{slanted}
  {\thspace}
  {\thspace}
  {\sl}
  {}
  {\bfseries}
  {.}
  {.5em}
  {}

\theoremstyle{slanted}

\newtheorem{prop}[equation]{Proposition}

\newtheorem{theo}[equation]{Theorem}
\newtheorem{lemm}[equation]{Lemma}
\newtheorem{coro}[equation]{Corollary}

\newtheoremstyle{roman}
  {\thspace}
  {\thspace}
  {\rm}
  {}
  {\it}
  {.}
  {.5em}
  {}

\theoremstyle{roman}

\newtheorem{defi}[equation]{Definition}
\newtheorem{rema}[equation]{Remark}

\newcommand{\daanmath}{\mathbb}
\newcommand{\zz}{\daanmath{Z}}

\newcommand{\rr}{\daanmath{R}}
\newcommand{\cc}{\daanmath{C}}
\newcommand{\qq}{\daanmath{Q}}

\newcommand\ra{\rightarrow}
\newcommand\longra{\longrightarrow}

\newcommand\sur{\mathrel{\to\kern-1.8ex\to}}
\newcommand\iso{\mathrel{\hookrightarrow\kern-1.8ex\to}}

\newcommand\longhookrightarrow{\lhook\joinrel\longrightarrow}

\newcommand\longsur{\mathrel{\longrightarrow\kern-1.8ex\to}}
\newcommand\longiso{\mathrel{\longhookrightarrow\kern-1.8ex\to}}

\newcommand\hugehookrightarrow{\lhook\joinrel\hugerightarrow}

\newcommand\hugesur{\mathrel{\hugerightarrow\kern-1.8ex\to}}
\newcommand\hugeiso{\mathrel{\hugehookrightarrow\kern-1.8ex\to}}

\newcommand{\diagrams}{\psset{arrowsize=2pt 4, arrowlength=.7, linewidth=.5pt, nodesep=4pt, arrows=->}}

\renewcommand{\>}{\rangle}

\newcommand{\barr}[1]{\overline{#1}}
\newcommand{\be}{\begin{equation}}
\newcommand{\ee}{\end{equation}}

\newcommand{\col}{\text{\upshape :\ }}

\begin{document} 
%


\title{The braid group of $\zz^n$}
\author{Daan Krammer}
\date{31 January 2007}
\address{Department of Mathematics, University of Warwick, Coventry CV4 7AL, UK}
\email{D.Krammer@warwick.ac.uk}
%
%

\begin{abstract} 
We define pseudo-Garside groups and prove a theorem about them parallel to Garside's result on the word problem for the usual braid groups. The main novelty is that the set of {\em simple\,} elements can be infinite. We introduce a group $B=B(\zz^n)$ which we call the braid group of $\zz^n$, and which bears some vague resemblance to mapping class groups. It is to $\GL(n,\zz)$ what the braid group is to the symmetric group $S_n$. We prove that $B$ is a pseudo-Garside group. We give a {\em small presentation\,} for $B(\zz^n)$ assuming one for $B(\zz^3)$ is given.
\end{abstract}


\maketitle


\thispagestyle{empty}

\begin{center} \parbox{.7\textwidth}{\tableofcontents} \end{center}


\section{Introduction}

Let $S$ be a compact oriented connected real $2$-manifold and $p$ a base point on the boundary of $S$. To keep things simple, let us define the {\em mapping class group\,} $M$ of $(S,p)$ as the group of automorphisms of $F:=\pi_1(S,p)$ coming from self-homeomorphisms of $S$ which fix $p$. Then $M$ acts on $F/F'=H_1(S,\zz)$. The kernel $I$ of this action is known as the Torelli group. We have an exact sequence
\[ 1\longra I \longra M \longra \Aut(F/F'). \]
In general, $M/I$ is infinite, and it is the symplectic group over the integers in the typical case where $S$ has just one boundary component.

If $S$ is a disk with $n$ holes then $M=B_n$, the braid group on $n$ strands. In this case, the Torelli group is also known as the pure braid group $P_n$. The quotient $B_n/P_n$ is finite (the symmetric group). Perhaps surprisingly, the pure braid group turns out to play a pivotal role in algebraically flavoured theories about $B_n$, for example Garside's greedy forms for braids \cite{gar69} and finite type invariants \cite{mw02}.

It would be interesting to generalise such theories to general mapping class groups $M$, see \cite{par05}. There are reasons why the role of the pure braid is expected to be taken by the Torelli group, especially Hain's infinitesimal presentation of the Torelli group \cite{hai97}. It seems hard to generalise Garside's theory to mapping class groups, which is why I propose to start at the other end. Which groups allow Garside type greedy forms and look a bit like mapping class groups?

Here is a geometric approach, which we don't pursue but may be helpful to think of. The braid group $B_n$ is the fundamental group of the space of $n$-element subsets of $\cc$. Let $A$ be the fundamental group of the space of additive subgroups of $\cc$ isomorphic to $\zz^n$. This group looks like the braid group: there are points moving around in the plane which aren't allowed to collide. We also have a surjection $A\ra\GL(n,\zz)$ which is similar to the surjection $M\ra M/I$.

Now $A$ seems less interesting. For one thing, it is huge and certainly not finitely generated. Which leads us to an algebraic approach.

The weak Bruhat ordering $<$ on the symmetric group $S_n$ is defined by $a\leq ab$ if and only if, for all $i,j\in\{1,\ldots,n\}$
\be (i<j\text{ and } iab<jab) \Rightarrow ia<ja. \label{zn108} \ee
The braid group $B_n$ can be presented by generators $\{r(a)\mid a\in S_n\}$ and relations $r(ab)=(ra)(rb)$ whenever $a\leq ab$.

Let $<$ be the standard lexicographic ordering on $\zz^n$. In analogy to (\ref{zn108}), define an ordering-like relation $\lesssim$ on $G:=\GL(n,\zz)$ by $a\lesssim ab$ if and only if, for all $x,y\in\zz^n$,
\[ (x<y\text{ and } xab<yab)\Rightarrow xa<ya. \]
We define the braid group of $\zz^n$, written $B=B(\zz^n)$, by generators $\{r(a)\mid a\in G\}$ and relations $r(1)=1$ and $r(ab)=(ra)(rb)$ whenever $a\lesssim ab$. Taken as a monoid presentation it yields the braid monoid $B^+$ of $\zz^n$.

The similarity between $B$ and the usual braid group $B_n$ is obvious. We have a surjection $B\ra\GL(n,\zz)$ which reminds us of $M\ra M/I$.

Our first main result is parallel to Garside's greedy form for braids and states that $B$ satisfies the conclusion of theorem~\ref{zn91}.

The braid group of $\zz^n$ is an example of what one may call {\em pseudo-Garside group\,} which is neither weaker nor stronger than what is called {\em Garside group\,} in \cite{deh02}. We define pseudo-Garside groups in definition~\ref{zn83}. Our second main result is that, again, Garside's theory can be generalised to pseudo-Garside groups (see theorem~\ref{zn91}). Of course, the paper deals with pseudo-Garside groups in general before it does the braid group of $\zz^n$.

There are two reasons why one needs different techniques for $B$ than for $B_n$. The first reason is that in fact, $\lesssim$ is not an ordering but what is known as a preordering. It turns out that this doesn't make the theory much different. The second and chief reason is that $G=\GL(n,\zz)$ is infinite, and indeed, has infinite chains. This makes it harder or impossible to use an approach based on a small presentation as is used in \cite{deh02} and other papers. Instead, we use the generators $r(a)$ from the beginning --- even in the definition of the braid group of $\zz^n$ as we saw. We need to build a theory of pseudo-Garside groups up from the ground which we do in section~\ref{zn97}.

Our third main result theorem~\ref{zn75} gives (in a precise sense) a presentation of $B^+(\zz^n)$ provided one knows one for $B^+(\zz^3)$. This result is similar to a result by Magnus \cite{mag34} which gives a presentation of $\GL(n,\zz)$ provided one knows one for $\GL(3,\zz)$. It is also analogous to the usual presentation (found by Artin) for the usual braid group.

It would be interesting to know if the braid group of $\zz^n$ has any use. Can the mapping class group be embedded in it?

In section~\ref{zn97} we study pseudo-Garside groups. In section~\ref{zn99} we prove enough results to conclude that $B$ (which is defined in section~\ref{zn102}) is an example of a pseudo-Garside group. In section~\ref{zn100} we obtain the small presentation for $B^+$.

\section{Lattices of total orderings} \label{zn99}

It is known that the weak Bruhat order (\ref{zn108}) on the symmetric group $S_n$ is a lattice ordering. In this section, we state and prove some analogous results. The main result of this section, and the only one needed in the sequel, is proposition~\ref{zn78}, and states that some ordering on the set of so-called lexicographic orderings on $\zz^n$ makes it into a lattice.

In the first subsection we make a lattice out of
the total orderings on a set. In the second subsection we specialise this by introducing a group action. In the third subsection we specialise even further and look at lexicographic orderings on $\zz^n$.

\subsection{The set-theoretic version}

Let $X$ be a set.
We write the set of total orderings on $X$ as
\[ \big\{ {\leq}_p\ \big|\ p\in P\big\} \]
where $P$ is an index set. We assume there is no repetition: ${\leq}_p\neq{\leq}_q$ whenever $p\neq q$. As usual each of these orderings $\leq_p$ comes with three more relations $\geq_p$, $<_p$ and $>_p$ whose meanings should be clear. We say that $p\in P$ has some property if $\leq_p$ has.

For $p\in P$ write $R_p=\{(x,y)\in X^2\mid x<_py\}$ (which equals ${<}_p$). Define
\begin{align*}
L_p\col P&\longrightarrow 2^{R_p}, \\
q &\longmapsto \big\{ (x,y)\in R_p\ \big|\ x>_q y\big\}=R_p\smallsetminus R_q.
\end{align*}
The image of $L_p$ is written $L_p(P)$. In this subsection we fix $p\in P$ and write ${<},R,L$ instead of ${<}_p,R_p,L_p$.

\begin{defi}
Call a set $A\subset R$ {\em closed\,} if for all $x,y,z\in X$ with $x<y<z$ one has
\[ (x,y)\in A\text{ and }(y,z)\in A\ \Longrightarrow\ (x,z)\in A. \]
Call it {\em co-closed\,} if $R\smallsetminus A$ is closed.
\end{defi}

\begin{lemm} \label{zn2} The map $L\col P\ra 2^R$ is injective and its image is the set of closed, co-closed subsets of $R$.
\end{lemm}

\begin{proof} Proof of injectivity of $L$. Let $q,r\in P$ be distinct. Then there are $x,y\in P$ with $x<_q y$ and $x>_r y$. We may assume $x<_p y$ (otherwise interchange $(x,q)$ with $(y,r)$). Then $(x,y)\not\in L(q)$ and $(x,y)\in L(r)$. This proves that $L$ is injective.

It is readily clear that $L(q)$ is closed and co-closed, for any $q$.

Let $A\subset R$ be closed and co-closed. We prove that $A\in L(P)$. Define a relation $<$ on $X$ by
\[ x<y \Leftrightarrow \big[
(x<_p y\text{ and } (x,y)\not\in A)
\text{ or }(y<_p x\text{ and } (y,x)\in A)
\big]. \]
A tedious case by case proof which we leave to the reader shows that $<$ is transitive. It follows readily that $<$ is an (anti-reflexive) total ordering. So ${<}={<}_q$ for some $q\in P$. Then $A=L(q)$ as required.
\end{proof}

\begin{defi} Define an ordering ${\leq}={\leq}^p$ on $P$ by
\be q\leq r\ \Longleftrightarrow\ \big(\text{for all }x,y\in X\col x<_py\text{ and }x<_ry\Longrightarrow x<_qy\big). \label{zn39} \ee
\end{defi}

For $p\in P$ we define $\barr{p}$ by ${\leq}_p={\geq}_{\barr{p}}$. It is clear that 
\be {\leq}^p={\geq}^{\barr{p}}. \label{zn18} \ee

\begin{lemm} \label{zn3} Let $p,q,r\in P$. Then $q\leq^p r$ if and only if $L_p(q)\subset L_p(r)$.
\end{lemm}

\begin{proof} Easy and left to the reader.
\end{proof}

A {\em lattice\,} is an ordered set such that any two elements $x,y$ have a least common upper bound or {\em join\,} $x\vee y$ and a greatest common lower bound or {\em meet\,} $x\wedge y$. A {\em complete lattice\,} is an ordered set such that any subset has a join and a meet.

\begin{prop} \label{zn4} Let $p\in P$. The ordered set $(P,\leq^p)$ is a complete lattice. For any subset $Q\subset P$, the set $L(\vee Q)$ is the closed subset of $R$ generated by $\cup_{q\in Q}L(q)$.
\end{prop}

\begin{proof} By lemmas~\ref{zn2} and \ref{zn3} we have an isomorphism of ordered sets $L\col P\ra L(P)$ where $L(P)$ is ordered by inclusion. We shall prove that $L(P)$ is a lattice.

By lemma~\ref{zn2}, $L(P)$ is the set of closed and co-closed subsets $A\subset R$. This is how we think of $L(P)$.

Let $M\subset L(P)$ be any subset. Let $B$ be the union of all elements of $M$ and let $C$ be the closure of $B$. Equivalently, for $(x,y)\in R$ we have $(x,y)\in C$ if and only if there exist $x=t_0< t_1< \cdots < t_n=y$ such that $(t_i,t_{i+1})\in B$ for all~$i$.

It remains to prove that $C$ is a join for $M$, because meets will follow through the symmetry ${\leq}^p={\geq}^{\barr{p}}$ in (\ref{zn18}). Even less is enough, namely, to prove that $C$ is co-closed.

Let $x<y<z$ ($x,y,z\in X$) and suppose $(x,z)\in C$. We want to prove $(x,y)\in C$ or $(y,z)\in C$. By construction there are $x=t_0<\cdots< t_n=z$ such that $(t_i,t_{i+1})\in B$ for all~$i$.

Suppose first $t_i=y$ for some $i$. Then $(x,y)=(t_0,t_i)\in C$.

Suppose next $t_i<y<t_{i+1}$ for some $i$. We know that $(t_i,t_{i+1})\in A$ for some $A\in M$. As $A$ is co-closed, it contains $(t_i,y)$ or $(y,t_{i+1})$, say, $(t_i,y)\in A$. Since $(x,t_i)=(t_0,t_i)\in C$ and $(t_i,y)\in A\subset C$ and $C$ is closed we conclude $(x,y)\in C$. The other case is similar and this proves that $C$ is co-closed as required.\end{proof}

By taking $X$ to be finite in proposition~\ref{zn4} one recovers that the weak Bruhat order (\ref{zn108}) on the symmetric group is a lattice.

\subsection{Group actions on $X$}

We retain the notation of the previous subsection, except that we won't assume any $p\in P$ to be fixed.

The following is obvious.

\begin{lemm} \label{zn43} Let $g$ be a permutation of $X$. If $p,q,r\in P$ are $g$-invariant then so are $q\vee_p r$ and $q\wedge_p r$.\qed
\end{lemm}

From now on we assume that $X=\zz^n$ where $n\geq 0$. An element $p\in P$ is said to be {\em translation invariant\,} if $x+z<_py+z$ $\Leftrightarrow$ $x<_p y$ for all $x,y,z\in \zz^n$.

\begin{lemm} \label{zn34} Let $p\in P$, $Q\subset P$. If $p$ and all elements of $Q$ are translation invariant then so are $\vee_p Q\in P$ and $\wedge_p Q\in P$ (which are defined by proposition~\ref{zn4}).
\end{lemm}

\begin{proof} Apply lemma~\ref{zn43}, letting $g$ range over all translations $\zz^n\ra\zz^n$, $x\mapsto x+y$ where $y\in\zz^n$.\end{proof}

\subsection{Lexicographic orderings}

We write $G=\GL(n,\zz)$ which acts on $\zz^n$ on the right.

\begin{defi}[Lexicographic] \label{zn79} Let $e_1$, \ldots, $e_n$ be the standard basis of $\zz^n$. We define the {\em standard lexicographical ordering\,} $\leq_\ell$ on $\zz^n$ as follows.
\[ y+\sum_{i=1}^n x_i\,e_i>_\ell y\ \Longleftrightarrow\  0=x_1=\cdots= x_{i-1}<x_i\text{ for some $i$.} \]
This ordering is total and translation invariant. We call $p\in P$ {\em lexicographic\,} if there exists $g\in G$ such that $x<_p y$ $\Leftrightarrow$ $xg<_\ell yg$ for all $x,y\in\zz^n$.
\end{defi}

\begin{lemm} \label{zn62} Let $p,q,r\in P$. If $p,q,r$ are lexicographic then so are $q\vee_p r$ and $q\wedge_p r$.
\end{lemm}

\begin{proof} See subsection~\ref{zn81}.
\end{proof}

Lemma~\ref{zn34} says that if $p\in P$ is translation invariant, then the set of translation invariant elements of $P$ is a complete sublattice of $(P,<^p)$. In particular, it is itself a complete lattice. Likewise, lemma~\ref{zn62} implies the following.

\begin{prop} \label{zn78} Let $p\in P$ be lexicographic. Then the set of lexicographic elements of $P$ has a lattice ordering $<^p$.\qed
\end{prop}

It is easy to show that the lattice of proposition~\ref{zn78} is not complete in general.

\subsection{Proof of lemma~\ref{zn62}} \label{zn81}

In this subsection we sketch a proof of lemma~\ref{zn62}. It can be skipped on first reading.

The {\em standard lexicographic ordering\,} ${<}_\ell$ and the {\em lexicographic orderings\,} on $\qq^n$ are defined just as for $\zz^n$ in definition~\ref{zn79}. Let $H_\qq\subset\GL(n,\qq)$ denote the group of linear automorphisms of $\qq^n$ preserving ${<}_\ell$; it is the group of upper triangular matrices in $\GL(n,\qq)$ with positive entries on the diagonal.

\begin{lemm} \label{zn109} (a). We have $\GL(n,\qq)=\GL(n,\zz)\cdot H_\qq$, that is, every element of $\GL(n,\qq)$ can be written $xy$ with $x\in\GL(n,\zz)$ and $y\in H_\qq$.

(b). There is a bijection from the set of lexicographic orderings on $\qq^n$ to those on $\zz^n$, defined by  ${<}\mapsto {<}\cap(\zz^n\times\zz^n)$.
\end{lemm}

\begin{proof} Part (b) is immediate from (a). 

Proof of (a). The inclusion $\supset$ is clear. We prove $\subset$ by induction on $n$. For $n=0$ there is nothing to prove. Assume it is true for $n-1$ and let $g\in\GL(n,\qq)$. The $\zz$-module generated by the entries of the first column of $g$ is of the form $a\zz$ ($a\in\qq_{>0}$). We may suppose that the first column of $g$ is zero, except that $g_{11}=a$ (if not, multiply $g$ on the left with a suitable element of $\GL(n,\zz)$). By the induction hypothesis there are $x,y$ such that $g=xy$, and $y\in H_\qq$, and $x\in\GL(n,\zz)$ preserves $e_1$ and $\zz e_2\oplus\cdots\oplus\zz e_n$. This finishes the proof of (a).
\end{proof}

Let $f\col\qq^n\ra\qq$ be $\qq$-linear and nonzero. Then $\{x\in\qq^n\mid f(x)>0\}$ and $\{x\in\qq^n\mid f(x)\geq0\}$ are called (respectively, open and closed) {\em half-spaces}. A {\em PL convex\,} set is an intersection of finitely many half-spaces (open or closed). Here PL stands for piecewise linear which should not be confused with piecewise affine. A subset of $\qq^n$ is said to be PL if it is a finite union of PL convex sets.

The following result is standard although I can't seem to find a reference.

\begin{prop} \label{zn63} Let $A,B\subset\qq^n$. If $A,B$ are PL then so are $\qq^n\smallsetminus A$, $A\cup B$ and $A+B:=\{a+b\mid a\in A,\ b\in B\}$.\qed
\end{prop}

A total translation invariant ordering $<$ on a $\qq$-vector space $V$ is called Archimedean if for all $x,y\in V$, if $x>0$ then $kx>y$ for some positive integer $k$. Equivalently, $(V,<)$ is isomorphic to a $\qq$-subspace of the real numbers with their standard ordering.

\begin{lemm} \label{zn64} Let $<$ denote a translation invariant total ordering on $\qq^n$. Then there exists a direct decomposition $\qq^n=V_1\oplus\cdots\oplus V_k$ and Archimedean orderings $<_i$ on $V_i$ such that the following holds. For all $v=\sum_{i=1}^k v_i$ ($v_i\in V_i$) one has $v>0$ if and only if there exists $i$ with
\be 0=v_1=\cdots=v_{i-1}<_i v_i. \label{zn66} \ee
\end{lemm}

\begin{proof} This is well-known but I can't seem to find a reference. It is also easily proved by the reader.
\end{proof}

For $p,q\in P$ translation invariant write 
\begin{align*}
K(p)&:=\{x\in\qq^n\mid x>_p 0\}, \\ 
U_p(q)&:=\{x\in\qq^n\mid x>_p0,\ x<_q0\}=K(p)\smallsetminus K(q), \\
U_p^0(q)&:=U_p(q)\cup\{0\}.
\end{align*}
From proposition~\ref{zn4} it follows that
\be U^0_p(q\vee_p r) = U^0_p(q)+U^0_p(r). \label{zn84} \ee

\begin{lemm} \label{zn65} Let $p\in P$ be translation invariant, and suppose that $K(p)$ is PL. Then $p$ is lexicographic.
\end{lemm}

\begin{proof} Let $<_0$ be the ordering on $\qq^n$ defined by $x<_0y$ if and only if $nx<_p ny$ for some integer $n>0$.

By the classification of translation invariant total orderings on $\qq^n$, lemma~\ref{zn64}, there is a direct decomposition $\qq^n=V_1\oplus\cdots\oplus V_k$ with $V_i\neq0$ and Archimedean orderings $<_i$ on $V_i$ such that the following holds. For all $v=\sum_{i=1}^k v_i$ ($v_i\in V_i$) one has $v>0$ if and only if there exists $i$ with (\ref{zn66}).

In order to prove the lemma suppose that, to the contrary, the ordering $<_p$ on $\zz^n$ is not lexicographic. By lemma~\ref{zn109}(b), $<_0$ isn't lexicographic either. Therefore, there exists $j$ such that $\dim V_j>1$. Choose a $2$-dimensional subspace $W\subset V_j$. Then $K_0(p):=K(p)\cap W$ is PL because $K(p)$ is. 

There exists a basis $(w_1,w_2)$ for $W$ and an irrational real number $\alpha\in\rr$ such that
\[ K_0(p)=\{xw_1+yw_2\mid x,y\in\qq,\ x+\alpha y>0\}. \]
On writing $\partial$ for the topological boundary for subsets of $W\otimes_\qq\rr$ and a bar for closures, it follows that
\[ \Big( \partial\,\overline{K_0(p)}\Big)\cap W=\{xw_1+yw_2\mid x,y\in\qq,\ x+\alpha y=0\}=\{0\}. \]
Since $K_0(p)$ is PL this implies $K_0(p)=\{0\}$ which is absurd.
\end{proof}

\paragraph{\bf Proof of lemma~\ref{zn62}.}  Let $p,q,r\in P$ be lexicographic and write $s=q\vee_p r$. By lemma~\ref{zn34}, $s$ is translation invariant.

Now $K(p)$ is easily seen to be PL (because $p$ is lexicographic). Similarly for $K(q)$ and $K(r)$. Moreover $U^0_p(q)=\{0\}\cup K(p)\smallsetminus K(q)$ is also PL, by proposition~\ref{zn63}. But
\[ U^0_p(s)=U^0_p(q)+U^0_p(r) \]
by (\ref{zn84}) which is again PL by proposition~\ref{zn63}. Now
\begin{align*} K(s)&=\{x\in\qq^n\mid x>_s0\} \\
&= \{x\in\qq^n\mid x>_s0,\ x>_p0\}\cup \{x\in\qq^n\mid x>_s0,\ x<_p0\} \\
&= \big[ K(p)\smallsetminus U(s) \big] \cup -U(s) \end{align*}
is PL. By lemma~\ref{zn65}, $s$ is lexicographic as required. Use the symmetry (\ref{zn18}) to deal with $q\wedge_p r$.\hfill$\Box$

\section{The braid group of $\zz^n$} \label{zn102}

In this section we introduce the braid group of $\zz^n$ and prove that it is, in the language of definition~\ref{zn83} below, a pseudo-Garside group.

\subsection{Notation and basics} \label{zn94}

\begin{defi}[Preorderings] \label{zn1} A {\em preordering\,} on a set $X$ is a relation $\lesssim$ satisfying transitivity ($x\lesssim y$ and $y\lesssim z$ imply $x\lesssim z$) and reflexivity ($x\lesssim x$ for all $x$). It follows that the relation $\sim$ defined by $x\sim y$ $\Leftrightarrow$ ($x\lesssim y$ and $x\gtrsim y$) is an equivalence relation, and the preordering induces an ordering $\leq$ on $X/{\sim}$ by $[x]\leq [y]$ $\Leftrightarrow$ $x\lesssim y$ where $[x]$ is the $\sim$-class of~$x$.
\end{defi}

Two totally ordered sets of the same cardinality are not necessarily isomorphic. However, $G=\GL(n,\zz)$ acts transitively on the set of {\em lexicographic\,} orderings on $\zz^n$. We'll gratefully make use of this fact which enables us to work with groups rather than groupoids. Groupoids are less convenient in notation though not by concept, and we could have dealt with groupoids had it been necessary.

We give $\zz^n$ the standard lexicographic ordering ${<}_\ell={<}$ (see definition~\ref{zn79}). Let $H$ denote the subgroup of $G$ of those elements preserving the ordering on~$\zz^n$.

Recall that the set of total orderings on $\zz^n$ is $\{ {\leq}_p\ \big|\ p\in P\}$. We define a map $u\col G\ra P$ by
\be 0<_{u(a)}x\Longleftrightarrow 0<xa \text{\quad for all $x\in\zz^n$}. \label{zn106} \ee

We define a relation $\lesssim$ on $G$ by
\be a\lesssim ab \quad \Longleftrightarrow \quad u(a)\leq u(ab) \label{zn38} \ee
where $\leq$ denotes $\leq^\ell$.
So by (\ref{zn39}), the definition of $\leq^\ell$,
\be a\lesssim ab \Longleftrightarrow \big( 0<x\text{ and }0<xab\Rightarrow 0<xa\text{ for all $x\in\zz^n$} \big). \label{zn40} \ee
From (\ref{zn38}) it is immediate that $\lesssim$ is a preordering on $G$.

As in definition~\ref{zn1} on preorderings, we define a relation $\sim$ on $G$ by $a\sim ab\Leftrightarrow (a\lesssim ab$ and $ab\lesssim a)$. So
\begin{align} a\sim ab &\Longleftrightarrow u(a)\leq u(ab)\leq u(a) \notag \\ &\Longleftrightarrow u(a)=u(ab)
\Longleftrightarrow
{<_{u(a)}} ={<_{u(ab)}}
\notag \\ & \Longleftrightarrow
\big(0<xa \Leftrightarrow 0<xab\text{ for all $x\in\zz^n$}\big) \notag
 \\ &\Longleftrightarrow
b\in H. \label{zn41}  \end{align}
So $G/{\sim}=G/H$. As in definition~\ref{zn1}, we have an ordering $\leq$ on $G/{\sim}=G/H$ given by $aH\leq bH$ $\Leftrightarrow$ $a\lesssim b$.

\begin{lemm} \label{zn90} The ordered set $(G/H,\leq)$ is a lattice with least element $H$ and greatest element $w_0H$ where $w_0=-1\in G$. For all $a,b\in G$ one has $a\lesssim b \Leftrightarrow w_0\,a\,w_0^{-1}\lesssim w_0\,b\,w_0^{-1}$.
\end{lemm}

\begin{proof} The bit involving $w_0$ is easy. The rest is just a reformulation of proposition~\ref{zn78}.
\end{proof}

\begin{defi} For $a,b\in G$ we write
\[ a*b=\left\{\begin{array}{@{}ll} ab & \text{if $a\lesssim ab$} \\[.5ex] \text{not defined\qquad} & \text{otherwise.} \end{array}\right. \]
\end{defi}

\begin{lemm} \label{zn61} Let $a,b,c\in G$. Then $(a*b)*c$ is defined if and only if $a*(b*c)$ is.
\end{lemm}

\begin{proof} See subsection~\ref{zn82}.
\end{proof}

In the language of definition~\ref{zn105} below, we have proved that $(G,H,\leq,w_0)$ is a pseudo-Garside germ. 

The definition of the {\em braid monoid $B^+$ of $\zz^n$} and the {\em braid group $B$ of $\zz^n$} is given in definition~\ref{zn83} below (it is put there because it can have a wider setting).

\subsection{Proof of lemma~\ref{zn61}} \label{zn82}

This subsection is devoted to a proof of lemma~\ref{zn61} and can be skipped in a first reading.

For $a\in G$ we write
\[ N(a)=\big\{ x\in\zz^n\smallsetminus\{0\}\ \big|\ 0<x\Leftrightarrow 0>xa\big\}. \]
For any two sets $A,B$ we write $A\oplus B$ for the set of elements in $A$ or $B$ but not both.

\begin{lemm} \label{zn59} For $a,b\in G$ we have $N(ab)=N(a)\oplus N(b)\,a^{-1}$.
\end{lemm}

\begin{proof} Let $x\in\zz^n$, $x>0$. Then
\begin{align*}
&x\in N(a)\oplus N(b)\,a^{-1} \\
\Longleftrightarrow{}&\big[ (0<x\Leftrightarrow 0>xa)\text{ or } (0<xa\Leftrightarrow 0>xab)\text{ but not both}\big] \\
\Longleftrightarrow{}&\big[ (0>xa)\text{ or } (0<xa\Leftrightarrow 0>xab)\text{ but not both}\big] \\
\Longleftrightarrow{}&0>xab
\Longleftrightarrow x\in N(ab).
\end{align*}
As $N(a)=-N(a)$ a similar result holds for negative $x$ and the proof is finished.
\end{proof}

\begin{lemm} \label{zn60} Let $a,b\in G$. Then the following are equivalent.
\begin{flalign*}
\quad  \text{\rm (1)\quad} & a\lesssim b. &\hfill \\
\text{\rm (2)\quad} & N(ab)=N(a)\cup N(b)\,a^{-1}. \\
\text{\rm (3)\quad} &  N(ab)=N(a)\sqcup N(b)\,a^{-1} \quad \text{(disjoint union).} \\
\text{\rm (4)\quad} & N(a)\cap N(b)\,a^{-1}=\varnothing.
\end{flalign*}
\end{lemm}

\begin{proof} By lemma~\ref{zn59}, (2), (3) and (4) are equivalent. The equivalence of (4) and (1) follows from
\begin{align*}
&N(a)\cap N(b)\,a^{-1}=\varnothing \\
\Longleftrightarrow{}& \nexists\, x\in\zz^n\col(0<x\Leftrightarrow 0>xa)\text{ and }(0<xa\Leftrightarrow 0>xab) \\
\Longleftrightarrow{}& \nexists\, x\in\zz^n\col 0<x\Leftrightarrow 0>xa\Leftrightarrow 0<xab \\
\Longleftrightarrow{}& \text{for all $x\in\zz^n$}\col (0<x\text{ and }0<xab)\Rightarrow 0<xa \\
\Longleftrightarrow{}& a\lesssim ab,
\end{align*}
the last equivalence using (\ref{zn40}).
\end{proof}

\paragraph{\bf Proof of lemma~\ref{zn61}.} We have
\begin{align*}
&(a*b)*c\text{ is defined} \Longleftrightarrow a\lesssim ab\text{ and }ab\lesssim abc \\
\stackrel{A}{\Longleftrightarrow}{}& N(ab) = N(a)\sqcup N(b)\,a^{-1} \text{ and } N(abc) = N(ab)\sqcup N(c)\,(ab)^{-1} \\
\stackrel{B}{\Longleftrightarrow}{}& N(abc)=N(a)\sqcup N(b)\,a^{-1} \sqcup N(c)\,(ab)^{-1} \\
\stackrel{B}{\Longleftrightarrow}{}& N(bc) = N(b)\sqcup N(c)\,b^{-1} \text{ and } N(abc) = N(a)\sqcup N(bc)\,a^{-1} \\
\stackrel{A}{\Longleftrightarrow}{}& b\lesssim bc\text{ and } a\lesssim abc \Longleftrightarrow a*(b*c)\text{ is defined}
\end{align*}
where $A$ indicates that we use (1) $\Leftrightarrow$ (3) in lemma~\ref{zn60} and $B$ that we use (4) $\Rightarrow$ (3).\hfill$\Box$

\section{Pseudo-Garside 
                 groups} \label{zn97}

\subsection{Summary}

\begin{defi} \label{zn105} A {\em pseudo-Garside germ\,} is a tuple $(G,H,{\leq},w_0)$ with the following properties.
\lista{{\rm(\PG\arabic{daana})}}{\setlength{\leftmargin}{4em}}
\item \label{zn32} Firstly, $H\subset G$ are groups, and $\leq$ is an ordering on $G/H$. We write $a\lesssim b$ $\Leftrightarrow$ $aH\leq bH$ ($a,b\in G$) and $a\sim b$ $\Leftrightarrow$ $aH=bH$. We require that $\leq$ is a lattice-ordering on $G/H$ with least element $H$ and greatest element $w_0H$. 
For $a,b\in G$ we write
\[ a*b=\left\{\begin{array}{@{}ll} ab & \text{if $a\lesssim ab$} \\[.5ex] \text{not defined\qquad} & \text{otherwise.} \end{array}\right. \]
We call $(x_1,\ldots,x_k)\in G^k$ {\em minimal\,} if $x_1*\cdots *x_k$ exists.%
\footnote{This is rather analogous to what \cite{bou68} calls {\em reduced decompositions\,} of elements of a Coxeter group.}
\item \label{zn33} Let $a,b,c\in G$. Then $(a*b)*c$ is defined if and only if $a*(b*c)$ is.
\item \label{zn77} For all $a,b\in G$ one has $a\lesssim b \Leftrightarrow w_0\,a\,w_0^{-1}\lesssim w_0\,b\,w_0^{-1}$.
\end{list}
\end{defi}

\begin{defi} \label{zn83} With a pseudo-Garside germ $(G,H,{\leq},w_0)$ we associate a monoid $B^+$ presented as follows.
\[ \begin{tabular}{ll}
Generators: & $\Omega=\{r(a)\mid a\in G\}$ (a copy of $G$). \\[.5ex]
Relations: & \parbox[t]{.6\textwidth}{$r(ab)=(ra)(rb)$ whenever $a,b\in G$ and $a*b$ is defined, that is, $a\lesssim ab$. Also, $r(1)=1$.} \end{tabular} \]
By $B$ we denote the group with the same presentation, taken as group presentation. We put $\Delta=rw_0$. We call $B^+$ a {\em pseudo-Garside monoid\,} and $B$ a {\em pseudo-Garside group}. Note that Garside groups in the sense of \cite{deh02} are not necessarily pseudo-Garside groups.
\end{defi}

One of the main results of this section is theorem~\ref{zn91} which says the following. In the above notation, every element of $B$ can be written $\Delta^kx_1\cdots x_\ell$ with $k\in\zz$, $\ell\geq 0$, $(x_1,\ldots,x_\ell)\in\Omega^\ell$ {\em strongly greedy\,} (see definition~\ref{zn93}) and $x_1\not\sim\Delta$ if $\ell>0$. Moreover, $k$ is unique and $(x_1,\ldots,x_\ell)$ is unique up to {\em strong equivalence\,} (see the beginning of subsection~\ref{zn92} for the missing definitions). This is very similar to one of Garside's results on the braid group \cite{gar69}.

In section~\ref{zn102} we proved that the braid group of $\zz^n$ is pseudo-Garside, so that it satisfies, for example, the conclusion of theorem~\ref{zn91}.

Of course, every group $G$ is pseudo-Garside: put $H=G$ so that also $B=G$. The challenge is to get $H$ small. In the case of the braid group of $\zz^n$, the group $H$ is nilpotent and therefore small for many purposes.

The remainder of this section is devoted to the proofs.

\subsection{Proofs} \label{zn92}

In this section we fix a pseudo-Garside germ $(G,H,{\leq},w_0)$ and we retain the notation of definitions~\ref{zn105} and~\ref{zn83}.

\begin{defi} 
Let $G^*$ be the free monoid on the set $G$. In order to keep the notation unambiguous, we identify $G^*$ with the disjoint union of Cartesian powers $\cup_{n\geq 0}G^n$. The unique element of $G^0$ is written $\varnothing$ or $()$. Elements of $G^1$ are often written $(a)$ rather than $a$ if $a\in G$. 

On $G^*$ we define a relation $\ra$ by
\[ {\ra}:=\left\{\big( u(a)(b*c)v,u(a*b)(c)v\big)\ \left|\ \begin{array}{@{}l@{}} u,v\in G^*,\ a,b,c\in G, \\ \text{$a*b$, $b*c$ defined} \end{array}\right\}.\right. \]
Let $\lesssim$ denote the reflexive-transitive closure of $\ra$. Clearly, $\lesssim$ is a preordering on $G^n$. Let $\sim$ denote the associated equivalence relation: $x\sim y$ $\Leftrightarrow$ $x\lesssim y\lesssim x$. Let $\approx$ denote the equivalence relation generated by $\ra$. In order to distinguish $\sim$ from $\approx$, we call 
$\approx$ the {\em equivalence\,} and
$\sim$ the {\em strong equivalence}.
\end{defi}

It is clear that $x\sim y$ $\Rightarrow$ $x\approx y$ ($x,y\in G^n$). One shouldn't confuse the preordering $\lesssim$ on $G^1$ (special case of $G^n$) with the preordering $\lesssim$ on $G$ as in (\PG\ref{zn32}).

Note that an element of $G^n$ is strongly equivalent to $(x_1,\ldots,x_n)$ if and only if it is of the form
\[ (x_1\,h_1^{-1}, h_1\,x_2\,h_2^{-1}, h_2\,x_3\,h_3^{-1}, \ldots, 
h_{n-2}\,x_{n-1}\,h_{n-1}^{-1}, h_{n-1}\,x_{n}) \]
for some $h_i\in H$.

Warning: If $x\ra y$ with $x\in G^m$ and $y\in G^n$ then $m=n$. The empty word $\varnothing\in G^0$ is {\em not\,} equivalent to $(1)\in G^1$. Only later will we identify the two.

\begin{lemm} \label{zn27} Let $a,b\in G$. Then $a\lesssim ab$ $\Leftrightarrow$ $b\lesssim a^{-1}w_0$.\qed
\end{lemm}

\begin{proof} For all $x\in G$, the expression $x*(x^{-1}w_0)$ is defined (and equals $w_0$) because $w_0$ is a greatest element. Therefore
\begin{align*}
a\lesssim ab
&\Longleftrightarrow
a*b\text{ is defined}
\Longleftrightarrow
(a*b)*(b^{-1}a^{-1}w_0)\text{ is defined}
\\ &\Longleftrightarrow
a*(b*(b^{-1}a^{-1}w_0))\text{ is defined}
\\ &\Longleftrightarrow
b*(b^{-1}a^{-1}w_0)\text{ is defined}
\Longleftrightarrow b\lesssim a^{-1}w_0.\qedhere \end{align*}
\end{proof}

\begin{lemm} \label{zn103} Let $a,b,c,d\in G$ be such that $bH\vee cH=dH$. If $a*b$ and $a*c$ are defined then so is $a*d$.
\end{lemm}

\begin{proof} We have $a\lesssim ab$ so, by lemma~\ref{zn27}, $b\lesssim a^{-1} w_0$. Similarly, $c\lesssim a^{-1} w_0$. Therefore
$d\lesssim a^{-1} w_0$. Using lemma~\ref{zn27} backwards yields $a\lesssim ad$ and therefore $a*d$ is defined.
\end{proof}

\begin{lemm} \label{zn24} Let $u,v,w\in G^*$, $u\ra v$, $u\ra w$. Then there exists $x\in G^*$ such that $v\ra x$, $w\ra x$.
\[ \psmatrix[rowsep=7mm, colsep=10mm]
u&v \\ w&x \endpsmatrix \diagrams \ncline{1,1}{1,2} \ncline{1,1}{2,1}
\ncline{2,1}{2,2}  \ncline{1,2}{2,2} \]
\end{lemm}

\begin{proof} First, consider the case 
$u=p\,u_0\,q$, $v=p\,v_0\,q$, $w=p\,w_0\,q$ where $u_0,v_0,w_0\in G^2$, $p,q\in G^*$. We may suppose $p=q=\varnothing$. Write $u=(a,b)$.

As $u\ra v$ we have
\[ u=(a,c*d)\ra (a*c,d)=v \]
for some $c,d\in G$. Likewise we can write
\[ u=(a,e*f)\ra (a*e,f)=w. \]
Write $r=c\vee e$, $r=c*p =e*q$, $b=r*s$. Since $r*s=(c*p)*s$ is defined, so is $p*s$. In fact, $cd=b=rs=cps$ so $d=p*s$. By lemma~\ref{zn103}, $a*r$ is defined, so $a*(c*p)$ is defined. We find
\[ v=(a*c,p*s)\ra (a*c*p,s)=(ar,s) \]
and likewise $w\ra (ar,s)$.

Next, consider the ``commutative'' case
\[ u=u_1\,u_2\,u_3\,u_4\,u_5,\quad
v=u_1\,v_0\,u_3\,u_4\,u_5,\quad
w=u_1\,u_2\,u_3\,w_0\,u_5 \]
where $u_2,u_4,v_0,w_0\in G^2$, $u_1,u_3,u_5\in G^*$. Then $x:=u_1\,v_0\,u_3\,w_0\,u_5$ does it.

It remains to consider the case
\begin{align*}
u &=p\,(a,u_2,u_3)\,q & (a,u_2,u_3)&\in G^3, \\
v &=p\,(v_1,c,u_3)\,q & (v_1,c,u_3)&\in G^3, \\
w &=p\,(a,w_2,e)\,q & (a,w_2,e)&\in G^3 \end{align*}
with $p,q\in G^*$. We may assume $p=q=\varnothing$. Since $u\ra v$ we have $u_2=b*c$, $v_1=a*b$ for some $b\in G$. As $u\ra w$ we have $u_3=d*e$, $w_2=u_2*d=(b*c)*d$ for some $d\in G$. In particular, $(b*c)*d$ is defined. By (\PG\ref{zn33}), $b*(c*d)$ is also defined. So we have the diagram
\[ \psmatrix[colsep=15mm, rowsep=6mm] \makebox[0mm][r]{$u={}$}(a,b*c,d*e) & (a*b,c,d*e)\makebox[0mm][l]{${}=v$} \\ \makebox[0mm][r]{$w={}$}(a,(b*c)*d,e) \\  (a,b*(c*d),e) & (a*b,c*d,e) \endpsmatrix \diagrams
\ncline{->}{1,1}{1,2}   \ncline{->}{3,1}{3,2}
\psset{nodesep=1pt} \ncline{->}{1,2}{3,2} \ncline{->}{1,1}{2,1} \ncline[doubleline=true, doublesep=1.8pt]{-}{2,1}{3,1}
\]
which finishes the proof.\end{proof}

\begin{defi} An element $x\in G^n$ is called {\em greedy\,} if its strong equivalence class is maximal, that is, $x\ra y$ implies $x\sim y$. We also say that $x$ is a {\em greedy form\,} of every element equivalent to it.
\end{defi}

Greedy elements (in an equivalence class) are not unique because all elements strongly equivalent to it are also greedy. But this is the only exception to uniqueness as we show now.

\begin{lemm} \label{zn26} (a). Every equivalence class $C\subset G^n$ has finite upper bounds, that is, for all $u,v\in C$ there exists $w\in C$ with $u\lesssim w$ and $v\lesssim w$.

(b). Every greedy element of $G^n$ is an upper bound (with respect to $\lesssim$) of all equivalent elements.
\end{lemm}

\begin{proof} (a). Let $u,v\in G^*$ be equivalent, that is, there exist
\[ u=u_0,\ u_1,\ \ldots,\ u_n=v, \qquad u_i\in G^* \]
such that for all $i$ one has $u_i\ra u_{i+1}$ or $u_{i+1}\ra u_i$. By induction on $n$, we prove that $\{u,v\}$ has an upper bound.

For $n=0$ there is nothing to prove. Assume it is true for $n-1$. Then $\{u_0,u_{n-1}\}$ has an upper bound $w$. If $u_n\ra u_{n-1}$ then $w$ is an upper bound of $u$ and $v$, so suppose $u_{n-1}\ra u_n$.

Since $u_{n-1}\lesssim w$ there exists a diagram as follows. 
\be \begin{aligned} \psmatrix[rowsep=7mm, colsep=9mm]
\makebox[0mm][r]{$u_{n-1}=:{}$}x_0 & x_1 & \cdots & x_k\makebox[0mm][l]{${}:=w$} \\ \makebox[0mm][r]{$v=u_{n}=:{}$}y_0 \endpsmatrix \diagrams
\ncline{1,1}{1,2} \ncline{1,2}{1,3} \ncline{1,3}{1,4} \ncline{1,4}{1,5} \ncline{1,1}{2,1} \end{aligned} \label{zn23} \ee
Using lemma~\ref{zn24} recursively, we can extend (\ref{zn23}) to a diagram as follows.
\[ \psmatrix[rowsep=7mm, colsep=9mm]
\makebox[0mm][r]{$u_{n-1}=:{}$}x_0 & x_1 & \cdots & x_k\makebox[0mm][l]{${}:=w$} \\ \makebox[0mm][r]{$v=u_{n}=:{}$}y_0 & y_1 & \cdots & y_k\makebox[0mm][l]{${}=:z$}  \endpsmatrix \diagrams
\ncline{1,1}{1,2} \ncline{1,2}{1,3} \ncline{1,3}{1,4} \ncline{1,4}{1,5} \ncline{2,1}{2,2} \ncline{2,2}{2,3} \ncline{2,3}{2,4} \ncline{2,4}{2,5} 
\ncline{1,1}{2,1} \ncline{1,2}{2,2} \ncline{1,4}{2,4} 
\]
So $v\lesssim z$ and also $u=u_0\lesssim w\lesssim z$.

(b). Immediate from (a).
\end{proof}

\begin{lemm} \label{zn29} Let $(1,x_1,\ldots,x_n)\leq(y_0,\ldots,y_n)$ (both in $G^{n+1}$) and suppose that $(x_1,\ldots,x_n)$ is greedy. Then $y_0\lesssim x_1$.
\end{lemm}

\begin{proof} The equivalence class $C$ of $(1,x_1,\ldots,x_n)$ contains $x:=(x_1,\ldots,x_n,1)$, which is greedy. By lemma~\ref{zn26}(b), $x$ is a greatest element in $C$. But $y:=(y_0,\ldots,y_n)$ is in $C$ too, so $y\lesssim x$. Therefore $y_0\lesssim x_1$.
\end{proof}

\begin{prop} \label{zn80} Every equivalence class in $G^*$ has a greedy element.
\end{prop}

\begin{proof} Let $A(n)$ denote the statement that every equivalence class in $G^n$ has a greedy element. We begin by proving $A(2)$. We tacitly use (\PG\ref{zn33}).

Let $(a,b)\in G^2$. Let $x\in G$ be such that $xH=a^{-1}w_0H\wedge bH$. There exists $y\in G$ with $b=x*y$. Now $x\lesssim a^{-1}w_0$ so by lemma~\ref{zn27}, $a*x$ is defined. So
\[ (a,b)=(a,x*y)\longra (a*x,y)=:t. \]
In order to prove that $t$ is greedy, suppose $t\ra t'$, say,
\[ t=(a*x,y)=(a*x,u*v)\longra ((a*x)*u,v)=t'. \]
We have $b=x*y=x*(u*v)=(x*u)*v$ whence
\be x*u\lesssim b. \label{zn28} \ee
As $a*(x*u)$ is defined we have $x*u\lesssim a^{-1}w_0$ by lemma~\ref{zn27} which we combine with (\ref{zn28}) to give
\[ a^{-1}w_0H\wedge bH=xH\leq (x*u)\,H \leq a^{-1}w_0H\wedge bH. \]
Therefore $u\sim 1$. This proves that $t$ is greedy and $A(2)$ is proved.

The proof of $A(n)$ is finished by induction on $n$. For $n\leq 1$ there is nothing to prove, and $A(2)$ has been proved above. We suppose $A(n-1)$ ($n\geq 3$) and aim to prove $A(n)$.

Let $C\subset G^n$ be an equivalence class. By $A(n-1)$, $C$ contains an element $x=(x_1,\ldots,x_n)$ such that $(x_2,\ldots,x_n)$ is greedy. By $A(2)$ there exists a greedy element $(y_1,y_2)\approx(x_1,x_2)$. By lemma~\ref{zn29} we have
\be \big[ x\ra(z_1,\ldots,z_n) \big] \Rightarrow z_1\lesssim y_1. \label{zn104} \ee
Define $w_1:=y_1$ and let $(w_2,\ldots,w_n)\approx(y_2,x_3,\ldots,x_n)$ be greedy. Note
\[ x\ra(y_1,y_2,x_3,\ldots,x_n)\ra w:=(w_1,\ldots,w_n). \]
In order to prove $w$ to be greedy, assume $w\ra z=(z_1,\ldots,z_n)$. By (\ref{zn104}) we have $z_1\sim w_1$ (because $y_1=w_1\lesssim z_1 \lesssim y_1$). By greediness of $(w_2,\ldots,w_n)$ we get $w\sim z$. So $w$ is greedy and the proof is finished.
\end{proof}

\begin{defi} \label{zn93} An element $(x_1,\ldots,x_n)\in G^n$ is called {\em strongly greedy\,} if it is greedy and ($n\geq 2$ $\Rightarrow$ $x_n\not\sim 1$) and ($n=1$ $\Rightarrow$ $x_1\neq 1$).
\end{defi}

We say that an element $(rx_1,\ldots,rx_k)\in\Omega^k$ (or two such) has some property (greedy, strongly greedy, equivalent, strongly equivalent) if $(x_1,\ldots,x_k)\in G^k$ has.

\begin{theo} \label{zn101} Every element of $B^+$ can be written $x_1\cdots x_n$ with $n\geq 0$, $x_i\in\Omega$ and $(x_1,\ldots,x_n)$ strongly greedy. Moreover, $(x_1,\ldots,x_n)$ is unique up to strong equivalence.
\end{theo}

\begin{proof} Let $\equiv$ denote the smallest equivalence relation on $G^*$ containing $\approx$ and such that $u(1)v\equiv uv$ for all $u,v\in G^*$. Then $B^+\cong G^*/{\equiv}$. We have
\[ x\sim y \Rightarrow x\approx y \Rightarrow x\equiv y \]
for all $x,y\in G^*$.

We shall define a map $S$ from $G^*$ to the set of strongly greedy elements in $G^*$. Let $R(x)$ denote any greedy element with $R(x)\approx x$ (it is not unique but we just choose one). Write $R(x)=(a_1,\ldots,a_n)$. If $n\leq 1$ we put $S(x):=R(x)$. If $n\geq 2$, let $k$ be maximal such that $a_k\not\in H$ and write $b=a_k\cdots a_n$. We put $S(x)=(a_1,\ldots,a_{k-1},b)$. Then $S(x)$ is strongly greedy and $S(x)\equiv x$. Also, if $x$ is strongly greedy then $x\sim R(x)=S(x)$. 

We claim that for all $x,y\in G^*$, if $x\equiv y$ then $S(x)\sim S(y)$. By the definition of $\equiv$, we need to prove this only if $x=u(1)v$, $y=uv$ (with $u,v\in G^*$) or if $x\approx y$. The case of $x\approx y$ is trivial. Now suppose $x=u(1)v$ and $y=uv$. Then $x\lesssim uv(1)=y(1)$ so $x\approx y(1)$ so $R(x)\sim R(y)(1)$ and $S(x)\sim S(y)$.

Now we can prove the lemma. Existence. Let $x\in G^*$. Then $S(x)$ is strongly greedy and $x\equiv S(x)$ as required.

Uniqueness. Let $x,y\in G^*$ be strongly greedy and $x\equiv y$. Because of $x\equiv y$ we get $S(x)\sim S(y)$. So $x\sim S(x)\sim S(y)\sim y$ as required.
\end{proof}

Note that we haven't used (\PG\ref{zn77}) so far. It is used in the proof of the following result.

\begin{theo} \label{zn91}
Every element of $B$ can be written $\Delta^kx_1\cdots x_\ell$ with $k\in\zz$, $\ell\geq 0$, $(x_1,\ldots,x_\ell)\in\Omega^n$ strongly greedy and $x_1\not\sim\Delta$ if $x_1$ is defined. Moreover, $k$ is unique and $(x_1,\ldots,x_\ell)$ is unique up to strong equivalence.
\end{theo}

\begin{proof} Easy using (\PG\ref{zn77}) and theorem~\ref{zn101} and left to the reader. \end{proof}

\section{A small presentation for $B^+$} \label{zn100}

The main result in this section is theorem~\ref{zn75} which gives a presentation of $B^+$, the braid monoid of $\zz^n$, in terms of generators and relations. Our approach is quite similar to Magnus' way \cite{mag34} to present $\GL(n,\zz)$ assuming that one has a presentation of $\GL(3,\zz)$.

\begin{defi} We define $G_i$ ($1\leq i<n$) to be the group of those $g\in G$ which preserve each $e_j$ ($j\not\in\{i,i+1\}$) as well as $\zz e_i\oplus\zz e_{i+1}$. Define $s_i\in G$ ($1\leq i\leq n$) by $e_i\,s_i=-e_i$ and $e_j\,s_i=e_j$ for all $j\neq i$. Note that $G_{i-1}\cap G_i=\<s_i\>$.
\end{defi}

\begin{lemm} \label{zn67} (a). Let $a\in G$, $b\in H$. Then $a*b$ and $b*a$ are defined.

(b). For all $a\in G_i$ and $x\in H$ there are $b\in G_i$, $y\in H$ such that
\be ax=yb. \label{zn68} \ee

(c). Same as (b) with $xa=by$ instead of (\ref{zn68}).
\end{lemm}

\begin{proof} (a). Clearly, $a*(b*b^{-1})$ and $(b^{-1}*b)*a$ are defined. By lemma~\ref{zn61}, $a*b$ and $b*a$ are defined.

Parts (b) and (c) are easy and left to the reader.
\end{proof}

\begin{figure}
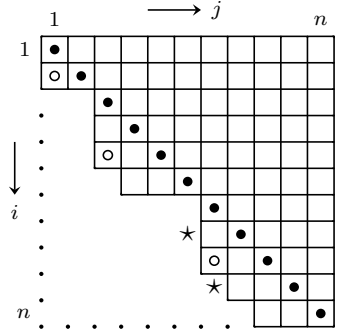
 \centering
\caption{A shape $C$. The black dots form the main diagonal. The white dots are the possible positions of $(i+1,j)\in C\smallsetminus A$ in proposition~\ref{zn57}. In case $(i+1,j)$ is the lower white dot, the stars are the positions $(i,j-1)$, $(i+2,j)$ which proposition~\ref{zn57} assumes not to be in $A$.\label{zn76}}
\[ \psset{linewidth=.6pt, unit=3.5mm, arrows=c-c, labelsep=4pt,
arrowsize=2pt 4, arrowlength=.7} 
\everypsbox{\scriptstyle}
\pspicture(-1,0)(11,12)$ \psline(0,11)(11,11) \psline(0,10)(11,10) \psline(0,9)(11,9) \psline(2,8)(11,8) \psline(2,7)(11,7) \psline(2,6)(11,6) \psline(3,5)(11,5) \psline(6,4)(11,4) \psline(6,3)(11,3) \psline(6,2)(11,2) \psline(7,1)(11,1) \psline(8,0)(11,0)
\multido{\ra=.5+1,\rb=10.5+-1}{11}{\pscircle*(\ra,\rb){1.8pt}}
\psline(0,9)(0,11) \psline(1,9)(1,11) \psline(2,6)(2,11) \psline(3,5)(3,11) \psline(4,5)(4,11) \psline(5,5)(5,11) \psline(6,2)(6,11) \psline(7,1)(7,11) \psline(8,0)(8,11) \psline(9,0)(9,11) \psline(10,0)(10,11) \psline(11,0)(11,11)
\pscircle(.5,9.5){2pt}
\pscircle(2.5,6.5){2pt}
\pscircle(6.5,2.5){2pt}
\uput[90]{0}(.5,11){1}
\uput[90]{0}(10.5,11){n}
\uput[180]{0}(0,.5){n}
\uput[180]{0}(0,10.5){1}
\multido{\ia=0+1}{9}{\pscircle*(0,\ia){.8pt}}
\multido{\ia=1+1}{7}{\pscircle*(\ia,0){.8pt}}
\rput(5.5,3.5){{\textstyle \star}} \rput(6.5,1.5){\textstyle \star}
\psline{->}(-1,7)(-1,5) \uput[-90]{0}(-1,5){i}
\psline{->}(4,12)(6,12) \uput[0]{0}(6,12){j}
$\endpspicture
\]
\end{figure}

\begin{defi}[shapes] See figure~\ref{zn76}. A {\em shape\,} is a set $A\subset\{1,\ldots,n\}^2$ such that $(i,i)\in A$ for all $i$, and for all $(i,j)\in\{1,\ldots,n\}^2$
\[ (i,j)\not\in A\quad\Longrightarrow\quad (i+1,j)\not\in A\text{ and }(i,j-1)\not\in A. \]
\end{defi}

Let $g\in G$. As usual, $g$ is a matrix $(g_{ij})_{ij}$ where $i,j$ range over $\{1,\ldots,n\}$; by definition (since $G$ acts on the right)
\[ e_i\,g=\sum_j g_{ij}\,e_j. \]
We define $\shape(g)$ (the shape of $g$) to be the smallest shape $A$ containing $\{(i,j)\mid g_{ij}\neq 0\}$. 

For a shape $A$, we define $G(A)$ to be the set of those elements of $G$ whose shape is contained in $A$. Note $H\,G(A)\,H=G(A)$. Note also $H\subset G(A_0)$ where $A_0$ is the smallest shape: $A_0=\{(i,j)\in\{1,\ldots,n\}^2\mid i\leq j\}$.

\begin{defi} \label{zn86} Recall $u\col G\ra P$ from (\ref{zn106}). For $g\in G$ we define $M(g)$ by
\begin{align*} M(g) &=\big\{x\in\zz^n\mid (0,x)\in L_\ell\big( u(g) \big) \big\}=\big\{ x\in\zz^n\mid x>0,\ x<_{ug}0\big\}
 \\ &= \big\{ x\in\zz^n\mid x>0>xg\big\}  \end{align*}
\end{defi}

\begin{lemm} \label{zn85} Let $a,b\in G$. Then $a\lesssim ab$ $\Leftrightarrow$ $M(a)\subset M(ab)$.
\end{lemm}

\begin{proof} We have
\begin{align*}
a\lesssim ab
&\stackrel{(\ref{zn40})}{\Longleftrightarrow} (\text{for all $x\in\zz^n$}\col 0<x\text{ and }0<xab\Rightarrow 0<xa) \\
&\Longleftrightarrow 
\big\{ x\in\zz^n\ \big|\ x>0>xa\big\} \subset
\big\{ x\in\zz^n\ \big|\ x>0>xab\big\} \\
&\Longleftrightarrow M(a)\subset M(ab). \qedhere
\end{align*}
\end{proof}

\begin{prop} \label{zn57} See figure~\ref{zn76}. Let $A\subset C$ be shapes with $\#C=\#A+1$. Suppose $(i+1,j)\in C\smallsetminus A$, $(i,j-1)\not\in A$, $(i+2,j)\not\in A$. Let $x\in G(C)$. 

(a). There are $y\in G_i$ and $z\in G(A)$ such that $x=y*z$.

(b). Any other pair $(y,z)$ with the same properties is of the form $(y,z)=(yt,t^{-1}z)$ with $t\in\<s_{i+1},H\cap G_i\>$.

\end{prop}

\begin{proof} (a). 
If $\shape(x)\subset A$ there is nothing to do (choose $y=1$, $z=x$), so suppose otherwise, that is, $x_{i+1,j}\neq 0$. Write 
\be
\begin{pmatrix} x_{ij} \\ x_{i+1,j} \end{pmatrix} =
\begin{pmatrix} au \\ bu \end{pmatrix} \label{zn45} \ee
where $a,b,u\in\zz$ with $a,b$ coprime and $u>0$. Choose $c,d\in\zz$ such that
\be ad-bc=1. \label{zn44} \ee
Define
\be y_0=\begin{pmatrix} av & cw \\ bv & dw \end{pmatrix} \label{zn74} \ee
where $v,w\in\{-1,1\}$ are to be determined later. Note that they are allowed to depend on $x$. Let $y\in G_i$ be the (unique) element of $G_i$ with $y_0$ in rows and columns of indices $i,i+1$. Put $z=y^{-1}x$. We need to show $x=y*z$ and $z\in G(A)$. 

We shall prove $z\in G(A)$. It is clear that $z=y^{-1}x\in G(C)$; we need to prove $z_{i+1,j}=0$. Consider the entries (\ref{zn45}) in $x$. The corresponding entries in $z$ are
\begin{align*}
y_0^{-1}\begin{pmatrix} x_{ij} \\ x_{i+1,j} \end{pmatrix}
&= y_0^{-1}\begin{pmatrix} au\\bu \end{pmatrix}
=u\,v\,w\begin{pmatrix} dw & -cw \\ -bv & av \end{pmatrix}\begin{pmatrix} a \\ b \end{pmatrix} \\ &=
u\,v\,w\begin{pmatrix} w(ad-bc) \\ 0 \end{pmatrix} = 
\begin{pmatrix} uv \\ 0 \end{pmatrix} \end{align*}
which shows that $z\in G(A)$.

It remains to show $x=y*z$, that is, $y\lesssim x$, or equivalently (by lemma~\ref{zn85}) $M(y)\subset M(x)$. Let $t\in M(y)$, that is (by definition~\ref{zn86}), $ty<0<t$. Write $t=\sum_k t_k\,e_k$. Note that $t_k=0$ for $k<i$ because otherwise, $t$ and $ty$ have the same sign, contradicting $ty<0<t$. For a similar reason 
\be (t_i,t_{i+1})\neq 0. \label{zn47} \ee

The $i$-th and $(i+1)$-th coefficients of $ty$ are
\begin{align*} \big((ty)_i,(ty)_{i+1}\big) &=(t_i,t_{i+1})\,y_0=(t_i,t_{i+1})\begin{pmatrix} av & cw \\ bv & dw \end{pmatrix} \\ &= \big(v\,(a\,t_i+b\,t_{i+1}),w\,(c\,t_i+d\,t_{i+1})\big) \end{align*}
which is a nonzero vector by (\ref{zn47}).

Recall that we have $t>0$ whose definition simplifies to
\be  t_i>0 \text{ or } (t_i=0,\ t_{i+1}>0). \label{zn46} \ee
Similarly we have $ty<0$ whose definition simplifies to
\be v\,(a\,t_i+b\,t_{i+1})<0 \text{ or } \big( a\,t_i+b\,t_{i+1}=0,\ w\,(c\,t_i+d\,t_{i+1})<0 \big). \label{zn51} \ee

We put $v=1$. Recall that $b\neq 0$; we assume $w\in\{-1,1\}$ has the sign of $-b$.

We shall prove
\be a\,t_i+b\,t_{i+1}\neq 0. \label{zn53} \ee
Suppose (\ref{zn53}) is false. By (\ref{zn51}) we find
\begin{align*}
0 &> b^2w\,(c\,t_i+d\,t_{i+1}) \\ &= b\,w\,\big( b\,c\,t_i+d\,(b\,t_{i+1})\big) \\
&= b\,w\,(b\,c\,t_i-a\,d\,t_i) \qquad \text{since not (\ref{zn53})} \\
&= -b\,w\,(a\,d-b\,c)\,t_i.
\end{align*}
But $-bw>0$ and $ad-bc=1$ so $t_i<0$, contradicting (\ref{zn46}). This proves (\ref{zn53}). By (\ref{zn51}) and (\ref{zn53}) we find
\be a\,t_i+b\,t_{i+1}<0. \label{zn54} \ee
The $i$-th coefficient of $tx$ is $at_i+bt_{i+1}<0$ by (\ref{zn54}) and all preceding coefficients are zero, so $tx<0$. Also $t>0$ so $t\in M(x)$ by definition~\ref{zn86}. This proves $M(y)\subset M(x)$ thus proving (a).

(b). We prefer to work with $y$ but not $z$; the conditions for $y$ to be satisfied are
\be y\lesssim x \label{zn72} \ee 
and 
\be y^{-1}x\in G(A). \label{zn73} \ee
An easy computation which we leave to the reader shows that $y$ satisfies (\ref{zn73}) if and only if there exists $p\in H\cap G_i$ such that $yp$ is of the form (\ref{zn74}) for some $v,w\in\{-1,1\}$ (or rather, its submatrix in rows and columns $i,i+1$). By lemma~\ref{zn67}, (\ref{zn72}) is also invariant under multiplying $y$ on the right with elements of $H\cap G_i$.

Assume therefore that $y$ is of the form (\ref{zn74}) and satisfies (\ref{zn72}). The proof will be finished by showing that $v=1$.

Suppose the contrary, $v=-1$. Choose $t_i,t_{i+1}\in\qq$ such that $t_i>0$ and $at_i+bt_{i+1}>0$, and put $t=t_ie_i+t_{i+1}e_{i+1}$. Then by (\ref{zn46}) we have $t>0$ and by (\ref{zn51}) we have $ty<0$. Therefore $t\in M(y)$. But the $i$-th coefficient of $tx$ is $at_i+bt_{i+1}>0$ by (\ref{zn54}) and all preceding coefficients are zero, so $tx>0$ and $t\not\in M(x)$. Therefore, $M(y)\not\subset M(x)$, that is, (\ref{zn72}) is false. This contradiction finishes the proof.
\end{proof}

\begin{lemm} \label{zn56} Let $x,y\in G$ be diagonal matrices (necessarily all diagonal entries being $1$ or $-1$). Suppose that for all $i$, if $e_ix=-e_i$ then $e_iy=-e_i$. Then $x\lesssim y$.
\end{lemm}

\begin{proof} By lemma~\ref{zn85} we need to prove $M(x)\subset M(y)$. Note $0\not\in M(x)$. Let $v=\sum_i v_i\,e_i\in\zz^n$, $v\neq 0$, say, $0=v_1=\cdots=v_{k-1}\neq v_k$. Then
\begin{align*}
v\in M(x) 
&\Leftrightarrow (v>0,\ vx<0)
\Leftrightarrow (v_k>0,\ e_k\,x=-e_k) 
\\ & \Rightarrow  (v_k>0,\ e_k\,y=-e_k)
\Leftrightarrow (v>0,\ vy<0)
\Leftrightarrow v\in M(y).
\qedhere \end{align*}
\end{proof}

\begin{coro} \label{zn107} The monoid $B^+$ is generated by
\be \bigg( \bigcup_{i=1}^{n-1} rG_i  \bigg) \cup rH. \label{zn55} \ee
\end{coro}

\begin{proof} Let $M\subset B^+$ denote the monoid generated by (\ref{zn55}). We know that $B^+$ is generated by $r(G)$ so we will be done if we prove that $r(x)\in M$ for all $x\in G$. We shall do this by induction on $\#\shape(x)$.

First suppose $\shape(x)$ is minimal, that is, $x$ is upper triangular. It is clear that
\[ x=t_1\cdots t_k\,h \]
for some $k\geq 0$, some distinct $t_1,\ldots,t_k\in\{s_1,\ldots,s_n\}$ and some $h\in H$. By lemma~\ref{zn56} it follows that $x=t_1*\cdots *t_k*h$ or equivalently
\[ rx=(rt_1)\cdots (rt_k)(rh). \]
But each $rt_i$ is in some $r(G_j)$, thus proving the statement if  $\shape(x)$ is minimal.

Assume $\shape(x)=C\neq A_0$ and assume that the result has been proved for all $z\in G$ with $\#\shape(z)<\#C$. The proof will be finished if we can prove the required result for $x$.

Note that there exist indices $i,j$ and a shape $A\subset C$ satisfying the assumptions of proposition~\ref{zn57}. For example, one can choose $j$ to be minimal such that the $j$-th column of $C$ differs from the $j$-th column of $A_0$; subject to this, let $i$ be maximal such that $(i,j)\in C$ and put $A=C\smallsetminus \{(i,j)\}$.

By proposition~\ref{zn57} there exist $y\in G_i$, $z\in G(A)$ such that $x=y*z$. Then $rx=(ry)(rz)$. Now $rz\in M$ by the induction hypothesis while $ry\in r(G_i)\subset M$ so $rx\in M$. The proof is finished.
\end{proof}

We define $H_i$ to be the group generated by $G_i$ and $H$. By lemma~\ref{zn67}, all its elements can be written $a*x$ and $y*b$ ($a,b\in G_i$, $x,y\in H$). We define $S$ to be the union of all $H_i$. We write $S^*=\cup_{n\geq 0} S^n$. A multiplication in $S^*$ is defined by concatenation, making it into a free monoid on $S^1$.

A {\em congruence\,} on a monoid $M$ is an equivalence relation $\sim$ on it such that the quotient set $M/{\sim}$ has a (necessarily unique) monoid structure such that the natural set map $M\ra M/\sim$ is a homomorphism of monoids.

Let $\sim$ denote the smallest congruence on $S^*$ satisfying the following.
\begin{enumerate}
\item[(S0)] We have $S^1\ni(1)\sim\varnothing\in S^0$.
\item[(S1)] We have $(x,y)\sim (xy)$ for all $x,y\in H_i$ such that $x*y$ is defined.
\item[(S2)] We have $(x,y)\sim (y,x)$ whenever $x\in G_i$, $y\in G_j$ and $|i-j|>1$.
\item[(S3)] We have $(x_1,y_1,x_2)\sim(y_2,x_3,y_3)$ whenever the following hold.
\begin{enumerate}
\item $x_k\in G_i$ for all $k$.
\item $y_k\in G_{i+1}$ for all $k$.
\item $x_1*y_1*x_2$ and $y_2*x_3*y_3$ are defined and equal.
\end{enumerate}
\end{enumerate}

\begin{lemm} \label{zn69} Consider the monoid homomorphism $f\col S^*\ra B^+$ defined by $f(x)=r(x)$ for all $x\in S^1=S$. For all $x,y\in S^*$, if $x\sim y$ then $f(x)= f(y)$.
\end{lemm}

\begin{proof} Let $x\in G_i$, $y\in G_j$, $|i-j|>1$. Then $xy=yx$. It is easy and left to the reader to prove that $x*y$ and $y*x$ are defined. So $f(x,y)=(rx)(ry)=r(xy)=r(yx)=(ry)(rx)=f(y,x)$. This proves that the map $f$ respects (S2). The other cases (S0), (S1), (S3) are trivial.
\end{proof}

Our aim is to prove the converse of lemma~\ref{zn69} ($f(x)=f(y)$ $\Rightarrow$ $x\sim y$) which we do in theorem~\ref{zn75}.

We write $T=\{1,2,\ldots,n-1\}$ and $T^*=\cup_{n\geq 0}T^n$ which, as $S^*$, is a free monoid on $T^1$ with concatenation as multiplication.

We say that $(x_1,\ldots,x_k)\in S^k$ has {\em type\,} $(y_1,\ldots,y_k)\in T^k$ if $x_i\in H_{y_i}$ for all $i$. Every word (= element in $S^*$) has at least one type, but possibly more.

Let $\ra$ be the smallest relation on $T^*$ with the following properties.
\begin{enumerate}
\item[(T0)] We have $T^0\ni\varnothing \ra (i)\in T^1$ for all $i$.
\item[(T1)] We have $(i,i)\ra(i)$ for all $i$.
\item[(T2)] We have $(i,j)\ra(j,i)$ whenever $|i-j|>1$.\
\item[(T3)] We have $(i,j,i)\ra(j,i,j)$ whenever $j=i+1$.
\item[(T4)] We have $axb\ra ayb$ whenever $x\ra y$ ($a,b,x,y\in T^*$).
\end{enumerate}
Notice the similarity with the congruence $\sim$ on $S^*$. If $t_1\ra t_2$ we say that $t_1$ can be {\em rewritten to\,} $t_2$.

The following lemma ties up $S$ with $T$.

\begin{lemm} \label{zn70} Let $t_1,t_2\in T^*$ with $t_1\ra t_2$. Then for every minimal word $w_1\in S^*$ of type $t_1$ there exists a minimal word $w_2$ of type $t_2$ such that $w_1\sim w_2$.
\end{lemm}

\begin{proof} It is enough to do this in the following cases.
\begin{enumerate}
\item[(0)] $T^0\ni\varnothing=t_1\ra t_2=(i)\in T^1$.
\item[(1)] $t_1=(i,i)\ra (i)=t_2$.
\item[(2)] $t_1=(i,j)\ra (j,i)=t_2$, $|i-j|>1$.
\item[(3)] $t_1=(i,j,i)\ra (j,i,j)=t_2$, $j=i+1$.
\end{enumerate}

Case (0). We have $w_1=\varnothing\in S^0$. Choose $w_2=1$. We have $w_1\sim w_2$ by (S0).

Case (1). Let $w_1\in S^*$ be minimal of type $(i,i)$. Write $w_1=(x,y)$. By minimality of $w_1$ then, $x*y$ is defined. So a good choice is $w_2=(xy)$ by (S1).

Case (2). Let $w_1$ be minimal and of type $(i,j)$, $|i-j|>1$. By lemma~\ref{zn67} we can write $w_1=(xa,by)$ where $x,y\in H$, $a\in G_i$, $b\in G_j$. Then $w_1=(xa,by)\sim(x,a,b,y)\sim (x,b,a,y)\sim(xb,ay)$ so $w_2=(xb,ay)$ is a good choice.

Case (3). Let $w_1\in S^*$ be minimal of type $(i,j,i)$ with $j=i+1$. By lemma~\ref{zn67} we can write $w_1=(a_1x_1,a_2x_2,a_3x_3)$. Using lemma~\ref{zn67} we can separate $G_k$ and $H$, that is, $w_1\sim w_3$ for some $w_3=(b_1,b_2,b_3,y)$ with $b_1,b_3\in G_i$, $b_2\in G_j$ and $y\in H$. By proposition~\ref{zn57} applied three times with $n=3$, we can write $b_1b_2b_3=c_1*c_2*c_3$ for some $c_2\in G_i$, $c_1,c_3\in G_j$. Then $w_3\sim w_2:=(c_1,c_2,c_3y)$ which is of type $(j,i,j)$ as required.
\end{proof}

\begin{rema} Note that case (3) in the proof of lemma~\ref{zn70} would fail if we replaced ``\,$t_1\ra t_2$'' by ``\,$t_2\ra t_1$'', or equivalently, interchanged $i$ and $j$ in the definition (T3) of $\ra$. Certainly, proposition~\ref{zn57} fails in that situation.
\end{rema}

\begin{lemm} \label{zn89} We have \[ (k,\ldots,1) (k,\ldots,1) \longrightarrow (k,\ldots,1) (k,\ldots,2) \] whenever $0< k<n$. (Note that, for example, $(i,j)(k,\ell)\in T^4$ is exactly $(i,j,k,\ell)$; we are using brackets here to ease reading).
\end{lemm}

\begin{proof} Induction on $k$. If $k=1$ it reads $(1,1)\ra(1)$ and follows from (T1). Suppose it is true for $k-1$. 

In the first arrow labelled (T2) in the following, we push the last $k$ to the left as far as possible using only (T2). Similarly, in the last arrow labelled (T2) the last $k$ is pushed back to the right. We write IH for the induction hypothesis.
{\allowdisplaybreaks
\begin{align*} 
&(k,\ldots,1) (k,\ldots,1) \\
\xrightarrow{\text{(T2)}} \ &
(k,k-1,k)(k-2,\ldots,1)(k-1,\ldots,1) \\
\xrightarrow{\text{(T0)}} \ &
(k,k-1,k)(k-1,\ldots,1)(k-1,\ldots,1) \\
\xrightarrow{\text{\ IH\ }} \ &
 (k,k-1,k)(k-1,\ldots,1)(k-1,\ldots,2) \\[1ex]
={}\ &(k)(k-1,k,k-1)(k-2,\ldots,1)(k-1,\ldots,2) \\
\xrightarrow{\text{(T3)}} \ &
(k)(k,k-1,k)(k-2,\ldots,1)(k-1,\ldots,2) \\ 
\xrightarrow{\text{(T1)}} \ &
(k,k-1,k)(k-2,\ldots,1)(k-1,\ldots,2) \\ 
\xrightarrow{\text{(T2)}} \ &
(k,\ldots,1) (k,\ldots,2). \qedhere \end{align*}}
\end{proof}

For fixed $n$, we write
\[ D_i=(n-1,\ldots,i)(n-1,\ldots,i+1)\cdots(n-1,n-2)(n-1) \in T^*. \]

\begin{prop} \label{zn87} Any element of $T^*$ can be rewritten to
\[ D_1=(n-1,\ldots,1) (n-1,\ldots,2)\cdots (n-1). \]
\end{prop}

\begin{proof} We use a double induction. Call the statement of the lemma $A(n)$. We prove $A(n)$ by induction on $n$. We clearly have $A(1)$. Assuming $A(n-1)$, we shall prove $A(n)$. Let $x\in T^k$. We prove $x\ra D_1$ by induction on $k$. 

For $k=0$ this follows from (T0). Assuming it true for $k-1$ we prove it for $k$. Since it is true for $k-1$ we can write $x=D_1(i)$ (product of $D_1$ and $(i)\in T^1$). If $i\neq 1$ then $A(n-1)$ tells us that $x=(n-1,\ldots,1)\,D_{2}(i)$ can be rewritten to $x=(n-1,\ldots,1)\,D_{2}=D_{1}$ as required. Suppose now $i=1$. Pushing the last letter $i=1$ as far as possible to the left using only (T2) we get
\[ \big[(n-1,\ldots,1)(n-1,\ldots,1)\big](n-1,\ldots,3)(n-1,\ldots,4)\cdots(n-1). \]
Rewriting the part in square brackets using lemma~\ref{zn89} yields $D_1$ as required.
\end{proof}

\begin{coro} Every element of $G$ is of the form $x_1\cdots x_k$ where $(x_1,\ldots,x_k)\in S^k$ is minimal of type $D_1$.
\end{coro}

\begin{proof} Immediate from corollary~\ref{zn107} and lemmas~\ref{zn70} and \ref{zn87}.
\end{proof}

\begin{lemm} \label{zn88} Let $a\in G_i$. Then $a$ and $as_{i+1}$ are comparable, that is, $a\lesssim as_{i+1}$ or $as_{i+1}\lesssim a$.
\end{lemm}

\begin{proof} Easy and left to the reader.
\end{proof}

\begin{theo} \label{zn75} (a). The converse to lemma~\ref{zn69} holds. In other words, the monoid $B^+$ is presented by generators $S$ and relations (S0)--(S3).

(b). The group $B$ is presented by the same generators and relations, taken as a group presentation.
\end{theo}

\begin{proof} Part (b) follows immediately from (a). We prove (a) by a double induction. Let $A(n)$ denote the statement of the theorem. Then clearly $A(2)$. Assuming $A(n-1)$ we prove $A(n)$ ($n>2$).

By the definition of $B^+$, it suffices to show that for {\em minimal\,} words $w_1,w_2\in S^*$, if $f(w_1)=f(w_2)$ then $w_1\sim w_2$. In order to keep notation simple, we repeatedly replace $w_1,w_2$ by equivalent minimal words until they are equal or obviously equivalent.

Let $X\col S^*\ra G$ be the natural homomorphism: $X(a)=a$ for all $a\in S^1=S$. Write $x=(x_{ij}):=X(w_1)=X(w_2)\in G$.

Proposition~\ref{zn87} tells us that any type of $w_1$ can be rewritten to standard type (that is, type $D_1$). By lemma~\ref{zn70}, $w_1$ is equivalent to a word of standard type. So we may now assume that $w_1$ is of standard type, say,
\[ w_1=(y_{n-1},\ldots,y_1)\,u \]
with $y_i\in H_i$ and $u\in S^*$ of type $D_2$. By lemma~\ref{zn67} we may even assume $y_i\in G_i$ (we can collect the necessary $H$ factors in the type $D_2$ factor $u$).

Here and henceforth we write $J(w)$ for the greatest $j$ such that the $(j,1)$-entry in $X(w)$ is nonzero ($w\in S^*$). Put $j=J(w_1)$.

If $n>k\geq j$ then $y_ky_{k-1}\cdots y_1\,X(u)$ has the same first column as $x$. Now $y_{j-1}\cdots y_1\,X(u)$ may have a different first column, but it cannot have a nonzero entry in position $(j+1,1)$, because that couldn't be cleaned up by removing any number of $y_i$ factors on the left. We may now assume $y_j\in\<s_{j+1},s_j\>$ (otherwise, move some $H$ factor into the type $D_2$ factor). Write $y_j=s_{j+1}^ps_j^q$ with $p,q\in\{0,1\}$. We may assume $y_j=s_j^q\in\<s_j\>$ because otherwise we replace $(y_{j+1},y_j)=(y_{j+1},s_{j+1}^ps_j^q)$ by $(y_{j+1}s_{j+1}^p,s_j^q)$. Then $y_i$ and $y_k$ commute whenever $i\leq j<k$. We may now assume $y_k=1$ for all $k>j$ (otherwise push them to the right into the type $D_2$ factor). Summarising, we now have $y_k=1$ for $k>j$ and $y_j\in\<s_j\>$.

We continue the proof by induction on $j$. First consider the case $j=1$. Then $x_{11}=(-1)^q$ for some $q\in\{0,1\}$. Moreover, $w_1=(s_1^q)u$, for some $u\in S^*$ of type $D_2$. Likewise, we have $w_2=(s_1^q)v$ for some $v\in S^*$ of type $D_2$. By $A(n-1)$, we have $u\sim v$ and therefore $w_1\sim w_2$. This establishes the case where $j=1$.

Supposing now the result for $j-1$ and smaller, we prove it for $j$.

Recall that $y_j\in\<s_j\>\subset G_{j-1}$. We may assume that $y_j=1$ because otherwise we replace $(y_j,y_{j-1})$ by $(1,y_jy_{j-1})$. Summarising, we have $w_1=(y_{j-1},\ldots,y_1)\,u$. Similarly, we can write $w_2=(z_{j-1},\ldots,z_1)\,v$. 

Let $C$ be the shape of $x=X(w_1)=X(w_2)$ and $A=C\smallsetminus \{(j,i)\}$ where $i=1$. Then $y=y_{j-1}$ and $y=z_{j-1}$ both satisfy the conditions of proposition~\ref{zn57}, that is, (\ref{zn72}) and (\ref{zn73}). By proposition~\ref{zn57}(b) we must have $z_{j-1}=y_{j-1}s_j^t$ for some $t\in\{0,1\}$ after moving some $H$ factors around.

We will now show that we may in fact suppose $t=0$. If not, then $z_{j-1}=y_{j-1}s_j$. By lemma~\ref{zn88}, $z_{j-1}$ and $y_{j-1}$ are comparable. After interchanging $w_1$ and $w_2$ if necessary, we may assume $y_{j-1}\lesssim z_{j-1}$, that is, $z_{j-1}=y_{j-1}*s_j$. In our word $w_2$, replace $z_{j-1}$ by $(y_{j-1},s_j)$. Now $s_j$ commutes with everything on its right but not in the type $D_2$ factor. Push $s_j$ into the type $D_2$ factor using (T2). Now $w_2$ is of the form as before except that $z_{j-1}=y_{j-1}$, that is, $t=0$.

So $w_1=(y_{j-1})v_1$ and $w_2=(y_{j-1})v_2$ where $J(v_1)$ and $J(v_2)$ are smaller than $j=J(w_1)=J(w_2)$. By the induction hypothesis we have $v_1\sim v_2$ and therefore $w_1\sim  w_2$. This proves (a).
\end{proof}



\begin{thebibliography}{B}

\bibitem[Bou68]{bou68} Bourbaki, Nicolas.
{\em Lie groups and Lie algebras. Chapters 4--6.}
Translated from the 1968 French original. Springer-Verlag, Berlin, 2002.

\bibitem[Deh02]{deh02} Dehornoy, Patrick. Groupes de Garside.  Ann.\ Sci.\ \'Ecole Norm.\ Sup.\ (4)  {\bf 35}  (2002),  no.\ 2, 267--306.

\bibitem[Gar69]{gar69} Garside, F.\ A.  The braid group and other groups. Quart.\ J.\ Math.\ Oxford Ser.\ (2)  {\bf 20}  (1969), 235--254.

\bibitem[Hai97]{hai97} Hain, Richard. Infinitesimal presentations of the Torelli groups. J.\ Amer.\ Math.\ Soc.\  {\bf 10}  (1997), no.\ 3, 597--651.

\bibitem[Mag34]{mag34} Magnus, Wilhelm. \"Uber $n$-dimensionale Gittertransformationen. Acta Math.\ {\bf 64} (1934), 353--367.

\bibitem[MW02]{mw02} Mostovoy, Jacob; Willerton, Simon. Free groups and finite-type invariants of pure braids. Math.\ Proc.\ Cambridge Philos.\ Soc.\ {\bf 132} (2002), no.\ 1, 117--130.

\bibitem[Par05]{par05} Paris, Luis. From braid groups to mapping class groups. \hfill\\ {\tt http://arxiv.org/abs/math.GR/0412024} .


\end{thebibliography}
\end{document}